\documentclass[11pt]{article}
\usepackage{amsfonts}
\usepackage{mathrsfs}
\usepackage{amsmath}
\usepackage{amssymb}
\usepackage{graphicx}
\newtheorem{theorem}{\scshape \mdseries  Theorem}[section]

\newtheorem{coro}[theorem]{\scshape \mdseries  Corollary}

\begin{document}

\title{\sf On Wiener polarity index of cactus graphs\thanks{
Supported by National Natural Science Foundation of China (11071002, 11126178), Program for New Century Excellent
Talents in University, Key Project of Chinese Ministry of Education
(210091), Specialized Research Fund for the Doctoral Program of
Higher Education (20103401110002), Science and Technological Fund of
Anhui Province for Outstanding Youth (10040606Y33), Project of Educational Department of Anhui Province (KJ2011A019),
Scientific
Research Fund for Fostering Distinguished Young Scholars of Anhui
University(KJJQ1001), Academic Innovation Team of Anhui University
Project (KJTD001B).}}
\author{Nan Chen, Wen-Xue Du, Yi-Zheng Fan\thanks{Corresponding
author.
 E-mail address: fanyz@ahu.edu.cn(Y.-Z. Fan), wenxuedu@gmail.com (W.-X. Du), 707666559@qq.com (N. Chen) }\\
  {\small  \it School of Mathematical Sciences, Anhui University, Hefei 230601, P. R. China} \\
 }
\date{}
\maketitle

\noindent {\bf Abstract:}
The Wiener polarity index of a graph $G$ is the number of unordered pairs of vertices $u, v$
such that the distance between $u$ and $v$ is $3$.
In this paper we give an explicit formula for the Wiener polarity index of cactus graphs.
We also deduce formulas for some special cactus graphs.

\noindent {\bf MR Subject Classifications:} 05C12, 92E10


\section{Introduction}
We use Trinajsti\'c \cite{Tri20} for terminology and notation.
Let $G=(V(G),E(G))$ be a connected graph.
The {\it distance} between two vertices $u$ and $v$ in $G$, denoted by $d_{G}(u,v)$, is the length of a shortest path between $u$ and $v$ in $G$.
The {\it Wiener polarity} index of a graph $G=(V,E)$, denoted by $W_{p}(G)$, is defined by $$ W_{p}(G):= \# \{\{u,v\}|d_{G}(u,v)=3,u,v\in V\}, \eqno(1)$$
which is the number of unordered pairs of vertices $\{u,v\}$ of $G$ such that $d_{G}(u,v)=3$.
In organic compounds, say paraffin, this number is the number of pairs of carbon atoms which are separated by three carbon-carbon bonds.
The name ``Wiener polarity index'' for the quantity defined in (1) is introduced by Harold Wiener \cite{Wie13} in 1947.
Wiener himself conceived the index only for acyclic molecules and defined it in a slightly different-yet equivalent-manner.
In the same paper, Wiener also introduced another index for acyclic molecules, called {\it Wiener index} or {\it Wiener distance index} and defined by
$$ W(G):=\sum_{\{u,v\} \subseteq V}d_{G}(u,v).$$
Wiener \cite{Wie13} used a linear formula of $W$ and $W_{p}$ to calculate the boiling points $t_{B}$ of the paraffins, i.e., $$t_{B}=aW+bW_{p}+c,$$
 where $a,b$ and $c$ are constants for a given isomeric group.

The Wiener index $W$ is popular in chemical literatures.
In the mathematical literature, it seems to be studied firstly by Entringer et al. \cite{Ent8} in 1976. From then on, many researchers studied the Wiener index in different ways. For instance, one can see \cite{Bar1}, \cite{Bon2}, \cite{Bon3}, \cite{Dan5}, \cite{Ent7}, \cite{Ent8}, \cite{Gut10}, \cite{Moh11} and \cite{Wie13} for the theoretical aspects, and \cite{Can4}, \cite{Gut9} and \cite{Moh12} for algorithmic and computational aspects.
Recently, Dobrynin et al. gave a comprehensive survey \cite{Dob6} for the Wiener index. The reader is referred to the paper for further details.

In the best of our knowledge, Wiener had some information about the applicability of this topological index.
Using the Wiener polarity index, Lukovits and Linert \cite{Luk14} demonstrated quantitative structure property relationships in a series of acyclic and cycle-containing hydrocarbons. Hosoya \cite{Hos15} found a physico-chemical interpretation of $W_{p}(G)$.
Recently, Du et al. \cite{Du16} described a linear time algorithm for computing the Wiener polarity index of trees and characterized the trees maximizing the index among all the trees of the given order.
Deng et al. \cite{Den17} characterized the extremal trees with respect to this index among all trees of order $n$ and diameter $d$. Deng \cite{Den18} also gave the extremal Wiener polarity index of all chemical trees with order $n$. Deng and Xiao \cite{Den19} found the maximum Wiener polarity index of chemical trees with $n$ vertices and $k$ pendants.

However, it seems that less attention has paid for Wiener polarity index of cycle-containing graph up to now.
While we are preparing this paper, we find that Behmaram et al. \cite{beh} discuss Wiener polarity index of fullerenes and hexagonal systems which contain no triangles or quadrangles,
and Hou et al. \cite{hou} discuss the maximum Wiener polarity index of unicyclic graphs.
In the paper we consider the Wiener polarity index of cactus graphs which are allowed  to have triangles or quadrangles or many cycles.

\section{Main result}

In this section, we introduce some graphs used in this paper.
Firstly, we introduce two graphs $G_{1}$ and $G_{2}$ as follows; see Fig. 2.1.

\begin{center}
\includegraphics[scale=.6]{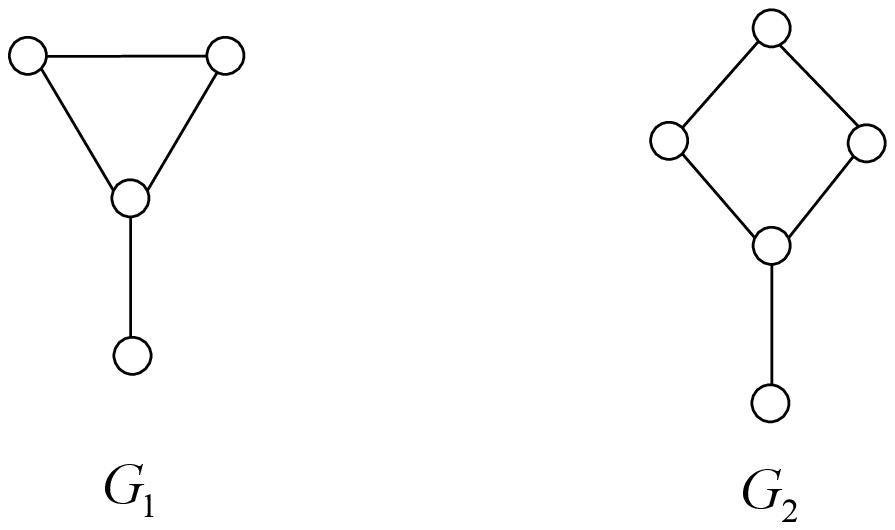}

\vspace{2mm}
{\small Fig. 2.1  The graphs $G_1$ and $G_2$ }
\end{center}

\noindent
Let $G$ be a graph. Denote by $c_{i}(G)$ the number of cycles of $G$ with length $i$.
The numbers of the induced subgraph of $G_{1}$ and $G_{2}$ in $G$ are denoted by $b_{1}(G)$ and $b_{2}(G)$, respectively.

A {\it cactus graph} $G$ is a connected graph in which no edge lies in more than one cycle.
 A {\it $k$-gon cactus graph} $G$ is a cactus graph in which every block is a $k$-gon or $C_k$, where $C_k$ denotes a cycle of length $k$.
 If each $k$-gon of a $k$-gon cactus $G$ has at most two cut-vertices, and  each cut-vertex is shared by exactly two hexagons,
 then $G$ is called a {\it chain $k$-gon cactus}.
 If, in addition, any two cut-vertices on a $k$-gon has distance $1$ (respectively, at least $2$),
 then this chain $k$-gon cactus is said {\it of type 1 } (respectively, {\it type 2}); see Fig. 2.2.
For a $k$-gon cactus of type 1 (respectively, type 2), expanding each of the cut-vertices to an edge, we will get
a graph called {\it ortho-chain $k$-gon cactus} (respectively, {\it meta-chain $k$-gon cactus}); see Fig. 2.3.

\begin{center}
\includegraphics[scale=.6]{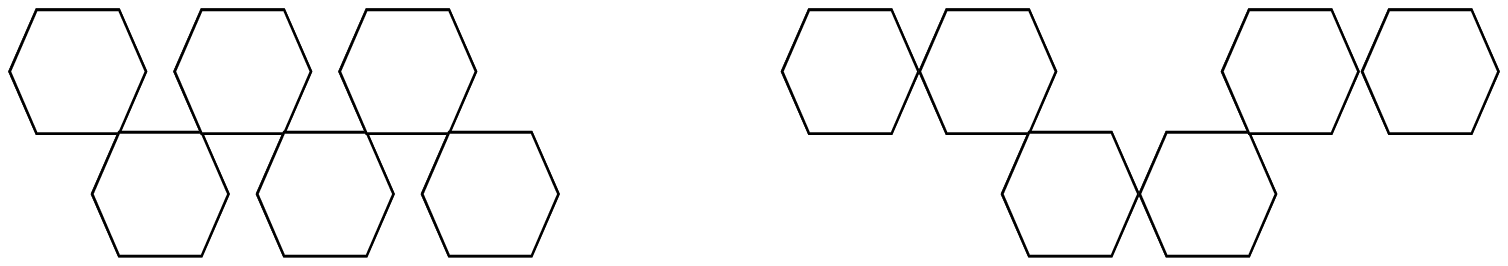}

\vspace{2mm}
{\small Fig. 2.2  Chain hexagonal cactuses of type 1 (left side) and type 2 (right side)}
\end{center}

\begin{center}
\includegraphics[scale=.6]{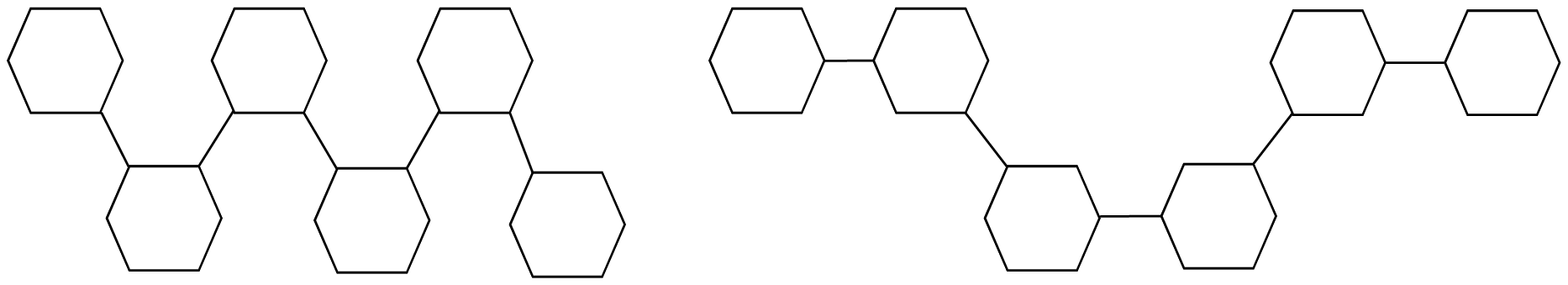}

\vspace{2mm}
{\small Fig. 2.3  Ortho-chain hexagonal cactus (left side) and meta-chain hexagonal cactus(right side)}
\end{center}

\begin{theorem} Let $G$ be a connected cactus graph. Then
$$W_{p}(G)=\sum_{uv\in E}(d_{G}(u)-1)(d_{G}(v)-1)-3c_{6}(G)-5c_{5}(G)-4c_{4}(G)-3c_{3}(G)-2b_{1}(G)-b_{2}(G).$$
\end{theorem}

{\bf Proof.}
We consider the edge $uv$ of the graph $G$ and choose a vertex $u'$ adjacent to $u$ and another vertex $v'$ adjacent to $v$.
Then we introduce a three-edge set in $G$, denoted by $Q$, defined by
$$Q:=\{(u',u,v,v')|u'u,uv,vv' \in E,u'\neq v,v' \neq u\}.$$
One can easily get $\#Q=\sum_{uv \in E}(d_{G}(u)-1)(d_{G}(v)-1)$.

If $d_{G}(u',v')=0$, namely $u',v'$ coincide and $u,u'=v',v$ in a common $C_{3}$, then $\#\{(u',u,v,v')\in Q| d_{G}(u',v')=0\}=3c_{3}(G)$.

If $d_{G}({u',v')}=1$, namely $u,u',v,v'$ lie in a common $C_{4}$, then $\#\{(u',u,v,v')\in Q|d_{G}(u',v')=1\}=4c_{4}(G)$.

If $d_{G}({u',v')}=2$, namely $u,u',v,v'$ lie in a common $C_{5}$, or $u,u',v $ lie in a common $C_{3}$ while $v'$ does not,
or $v,v',u$ lie in a common $C_{3}$ while $u'$ does not,
then $\#\{(u',u,v,v')\in Q|d_{G}(u',v')=2\}=5c_{5}(G)+2b_{1}(G)$.

If $d_{G}({u',v')}=3$, then $u,u',v,v' $ may lie in a common $C_{6}$,  or $u,u',v $ may lie in a common $C_{4}$ while $v'$ does not,
or $v,v',u $ may lie in a common $C_{4}$ while $u'$ does not,
then $$\# \{\{u',v'\}|d_{G}(u',v')=3,u,v \in V\}=\#\{(u',u,v,v')\in Q_{uv}|d_{G}(u',v')=3\}-3c_{6}(G)-b_{2}(G).$$

Combining the above discussion, we have
\begin{align*}
W_{p}(G)&=\#Q- \sum_{i=0}^2 \#\{(u',u,v,v')\in Q_{uv}| d_{G}(u',v')=i\}-3c_{6}(G)-b_{2}(G)\\
&=\sum_{uv\in E}(d_{G}(u)-1)(d_{G}(v)-1)-3c_{6}(G)-5c_{5}(G)-4c_{4}(G)-3c_{3}(G)-2b_{1}(G)-b_{2}(G).
\end{align*}
\hfill$\blacksquare$

\begin{coro} Let $G$ be a connected cactus graph such that every triangle or quadrangle has exactly one neighbor. Then
$$W_{p}(G)=\sum_{uv\in E}(d_{G}(u)-1)(d_{G}(v)-1)-3c_{6}(G)-5\sum_{i=3}^{5}c_{i}(G).$$
\end{coro}

{\bf Proof.}
Put $c_{3}(G)=b_{1}(G)$ and $c_{4}(G)=b_{2}(G)$ in Theorem 2.1.\hfill$\blacksquare$

\begin{coro}Let $G=(V,E)$ be a chain $k$-gon cactus of type 1 with $h \;(h \geq 2)$ $k$-gons. Then
$$
W_{p}(G) = \begin{cases}
4h-8, &\hbox{if~~} k=3,\\
8h-12, &\text{if~~} k=4,\\
12h-16, &\text{if~~} k=5,\\
15h-16, &\text{if~~} k=6,\\
(k+12)h-16, &\text{if~~} k \geq 7.
\end{cases}
$$
\end{coro}

{\bf Proof.}
It is easy to get $$\sum_{uv\in E}(d_{G}(u)-1)(d_{G}(v)-1)=(k+12)h-16.$$
So, by Theorem 2.1,
if $k=3$, then
$$W_{p}(G)=\sum_{uv\in E}(d_{G}(u)-1)(d_{G}(v)-1)-3c_{3}(G)-2b_{1}(G)=(15h-16)-3h-2(4h-4)=4h-8;$$
and if $k=4$, then
$$W_{p}(G)=\sum_{uv\in E}(d_{G}(u)-1)(d_{G}(v)-1)-4c_{4}(G)-b_{2}(G)=(16h-16)-4h-(4h-4)=8h-12.$$
Similarly, we get the remaining result and omit the details.\hfill$\blacksquare$

\begin{coro}Let $G=(V,E)$ be a chain $k$-gon cactus of type 2 with $h\;(h \geq 2)$ $k$-gons ($k \ge 4$). Then
$$
W_{p}(G) = \begin{cases}
4h-4, &\text{if~~} k=4,\\
8h-8, &\text{if~~} k=5,\\
11h-8, &\text{if~~} k=6,\\
(k+8)h-8, &\text{if~~} k \geq 7.
\end{cases}
$$
\end{coro}

{\bf Proof.}
One can get $$\sum_{uv\in E}(d_{G}(u)-1)(d_{G}(v)-1)=(k+8)h-8.$$
So,
if $k=4$, $W_{p}(G)=\sum_{uv\in E}(d_{G}(u)-1)(d_{G}(v)-1)-4c_{4}(G)-b_{2}(G)=(12h-8)-4h-(4h-4)=4h-4$.
The remaining proof is omitted.\hfill$\blacksquare$

\begin{coro}Let $G=(V,E)$ be an ortho-chain $k$-gon cactus with $h\;(h \geq 2)$ $k$-gons. Then
$$
W_{p}(G) = \begin{cases}
5h-6, &\text{if~~} k=3,\\
7h-8, &\text{if~~} k=4,\\
9h-10, &\text{if~~} k=5,\\
12h-10, &\text{if~~} k=6,\\
(k+9)h-10, &\text{if~~} k \geq 7.
\end{cases}
$$
\end{coro}

{\bf Proof.}
One can get $$\sum_{uv\in E}(d_{G}(u)-1)(d_{G}(v)-1)=(k+9)h-10.$$
So, if $k=3$, $W_{p}(G)=(12h-10)-3h-2(4h-4)=5h-6$; and
if $k=4$, $W_{p}(G)=(13h-10)-4h-(4h-4)=7h-8$.
The remaining proof is omitted.\hfill$\blacksquare$

\begin{coro}Let $G=(V,E)$ be a meta-chain $k$-gon cactus with $h\;(h \geq 2)$ $k$-gons. Then
$$
W_{p}(G) = \begin{cases}
6h-6, &\text{if~~} k=4,\\
8h-8, &\text{if~~} k=5,\\
11h-8, &\text{if~~} k=6,\\
(k+8)h-8, &\text{if~~} k \geq 7.
\end{cases}
$$
\end{coro}

{\bf Proof.}
One can get $$\sum_{uv\in E}(d_{G}(u)-1)(d_{G}(v)-1)=(k+8)h-8.$$
So, if $k=4$, $W_{p}(G)=(12h-8)-4h-(4h-4)=6h-6$.
The remaining proof is omitted.\hfill$\blacksquare$

\end{document}